\documentclass[11pt, twoside, a4paper]{article}
\usepackage[english]{babel}
\usepackage{amsmath,amssymb,amsthm,epsfig}
\newtheorem{theorem}{Theorem}
\newtheorem{lemma}[theorem]{Lemma}
\newtheorem{proposition}[theorem]{Proposition}

\newtheorem{definition}[theorem]{Definition}

\renewenvironment{proof}[1][.]{%
\bigskip\noindent{\bf Proof#1 }}{%
\hfill$\blacksquare$\bigskip}

\begin{document}

\title{ Entropy and variational principles for holonomic probabilities of IFS }
\author{A. O. Lopes (UFRGS) and Elismar R. Oliveira (UFRGS)}
\date{\today}
\maketitle

\begin{abstract}

An IFS ( iterated function system), $([0,1], \tau_{i})$, on the
interval $[0,1]$, is a family of continuous functions
$\tau_{0},\tau_{1}, ..., \tau_{d-1} : [0,1] \to [0,1]$.

Associated to a IFS one can consider a continuous map
$\hat{\sigma} : [0,1]\times \Sigma \to [0,1]\times \Sigma$,
defined by $\hat{\sigma}(x,w)=(\tau_{X_{1}(w)}(x), \sigma(w))$
were $\Sigma=\{0,1, ..., d-1\}^{\mathbb{N}}$, $\sigma: \Sigma \to
\Sigma$ is given by
$\sigma(w_{1},w_{2},w_{3},...)=(w_{2},w_{3},w_{4}...)$ and $X_{k}
: \Sigma \to \{0,1, ..., n-1\}$ is the projection on the
coordinate $k$.

A $\rho$-weighted system, $\rho \geq 0$, is a weighted system
$([0,1], \tau_{i}, u_{i})$ such that there exists a positive
bounded function $h : [0,1] \to \mathbb{R}$ and a probability $\nu $
on $[0,1]$ satisfying $ P_{u}(h)=\rho h, \quad P_{u}^{*}(\nu)=\rho
\nu$.

There is no meaning to ask if the probabilities $\nu $  on $[0,1]$
arising in IFS are invariant for a dynamical system, but, we can
ask if  probabilities $\hat{\nu}$ on $[0,1]\times \Sigma$ are
holonomic for $\hat{\sigma}$.

A probability  $\hat{\nu}$ on $[0,1]\times \Sigma$ is called
holonomic for $\hat{\sigma}$ if $ \int\, g \circ \hat{\sigma}\,
d\hat{\nu}=  \int \,g \,d\hat{\nu}, \, \forall g \in C([0,1])$. We
denote the set of holonomic probabilities by ${\cal H}$.

Via disintegration, holonomic probabilities $\hat{\nu}$ on
$[0,1]\times \Sigma$  are naturally associated to a
$\rho$-weighted system. More precisely,  there exist  a
probability $\nu$ on $[0,1]$ and $u_i, i\in\{0, 1,2,..,d-1\}$ on
$[0,1]$, such that is $P_{u}^*(\nu)=\nu$.

We consider holonomic ergodic probabilities and present the corresponding Ergodic Theorem (which is just an adaptation of a previous result by J. Elton).

For a holonomic probability $\hat{\nu}$ on $[0,1]\times \Sigma$ we
define the entropy $h(\hat{\nu})=\inf_{f \in \mathbb{B}^{+}} \int
\ln(\frac{P_{\psi}f}{\psi f}) d\hat{\nu}\geq 0$, where, $\psi  \in
\mathbb{B}^{+}$ is a fixed (any one) positive potential.

Finally, we analyze the problem: given $\phi \in \mathbb{B}^{+}$,
find the solution of the maximization problem
$$p(\phi)=
\sup_{\hat{\nu} \in \mathcal{H}} \{ \,h(\hat{\nu}) + \int
\ln(\phi) d\hat{\nu} \,\}.$$

We show na example where such supremum is attained by a holonomic not-invariant probability.
\end{abstract}

\vspace {.8cm}

\newpage
\section{ IFS and holonomic probabilities}

We want to analyze, in the setting of holonomic probabilities
\cite{GL} associated to an IFS, the concepts of entropy and
pressure.

We point out that this is a different problem from the usual one among invariant probabilities (see remarks 3 and 4 in section 7).

The present work is part of the PhD thesis of the second
author \cite{O}.

Our main point of view is the following: the study of the holonomic probabilities allows one
to understand all the transference operators $P_u$ and the associated stationary states when the IFS is considered as a realization of a Stochastic Process.

\begin{definition}\label{IFS}
An IFS (iterated function system), $([0,1], \tau_{i})$, on the
interval $[0,1]$, is a family of continuous functions
$\tau_{0},\tau_{1}, ..., \tau_{d-1} : [0,1] \to [0,1]$.
\end{definition}

Associated to a IFS one can consider a continuous map
$\hat{\sigma} : [0,1]\times \Sigma \to [0,1]\times \Sigma$,
defined by
$$\hat{\sigma}(x,w)=(\tau_{X_{1}(w)}(x), \sigma(w)),$$
were $\Sigma=\{0,1, ..., d-1\}^{\mathbb{N}}$, $\sigma: \Sigma \to
\Sigma$ is given by
$\sigma(w_{1},w_{2},w_{3},...)=(w_{2},w_{3},w_{4}...)$ and $X_{k}
: \Sigma \to \{0,1, ..., n-1\}$ is the projection on the
coordinate $k$. In this way one can see such system as a
Stochastic Process \cite{CR} \cite{BJ} \cite{BarnDemElton}
\cite{Elton} \cite{MU} \cite{DutkayJorg} \cite{Ja1} \cite{Ja2}.

If we consider a IFS as a multiple dynamical systems (several
maps) then, for a single point $x$ there exists several
combinations of ``orbits" on the IFS (using different $\tau_i$).
Considering the map $\hat{\sigma}$ one can describe the global
behavior of iterates of $x$. Moreover, one can think the IFS, as a
branching process with index in $\Sigma$. More precisely, we define
the
$n$-branch from $x \in [0,1]$ by $w \in \Sigma$, as
$$Z_{n}(x,w)=\tau_{X_{n}(w)} \circ \tau_{X_{n-1}(w)} \circ ...
\tau_{X_{1}(w)} (x).$$ With this notation, we have
$$\hat{\sigma}^{n}(x,w)=(Z_{n}(x,w),\sigma^{n}(w)).$$

\begin{definition}\label{WeightSystem}
A weighted system (see \cite{Stenflo}, Pg. 6) is a triple,
$([0,1], \tau_{i}, v_{i})$, were $([0,1], \tau_{i})$ is a IFS
where the $v_{i}$'s, $i\in\{0,1,...,d-1\}$, are measurable and nonnegative bounded maps.
The condition $\sum_{i=0}^{d-1}v_{i}(x)=1$ is not required.
\end{definition}

In some examples the $u_i$, $i\in\{0,1,...,d-1\}$, come
from a measurable  bounded potential $\phi : [0,1] \to [0,+
\infty)$, that is,
$$u_{i}(x)=\phi(\tau_{i}(x)), \; \forall i=0,...,d-1\,.$$
The function $\phi$ is called \textit{weight function}, (in the
literature this function is also called $g$-function, see
\cite{Keane} for example). Note that $\phi$ can attain the value $0$. This is useful for some applications of IFS to wavelets \cite{DutJorgDisint} \cite{DutkayJorg}. We do not assume in this general definition that $\sum_{i=0}^{d-1}u_{i}(x)=1, \quad \forall x \in [0,1]$.

\begin{definition}\label{IFSwithProbabilities}
A IFS with probabilities, $([0,1], \tau_{i}, u_{i})$, $i\in\{0,1,...,d-1\}$, is a IFS
with a vector of measurable nonnegative bounded functions on
$[0,1]$,
$$u(x)=(u_{0}(x), u_{1}(x), ..., u_{d-1}(x)),$$ such that,
$\sum_{i=0}^{d-1}u_{i}(x)=1, \quad \forall x \in [0,1]$.
\end{definition}

\begin{definition}\label{UniformSystem}
A IFS with probabilities $([0,1], \tau_{i}, u_{i})$ is called
``uniform normalized" if there exists  a  weight function $\phi$ such that
$u_{i}(x)=\phi(\tau_{i}(x)), \; \forall i=0,...,d-1$ and
$$\sum_{i=0}^{d-1}\phi(\tau_{i}(x))=1.$$

In this case,
we write the IFS as $([0,1], \tau_{i}, u_{i})=([0,1], \tau_{i}, \phi)$.

\end{definition}

The above definition  is a strong restriction in the weighted system. Several problems in the
classical theory of Thermodynamic Formalism for the shift or for a
$d$ to $1$ continuous expanding transformations $T: S^1 \to S^1$ can
be analyzed via a IFS with a weight function $\phi$ (see \cite{Parry}). In
this case the $\tau_i, i\in\{0,2,..,d-1\}$, are the inverse branches
of $T$.

We will consider later the pressure problem for a weight function $\phi$ which is not necessarily uniformly normalized.

We now return to the general case.

\begin{definition}\label{TransfOperator}
Given a weighted system, $([0,1], \tau_{i}, u_{i})$, we will
define de Transference Operator (or Ruelle Operator)   $P_{u}$ by
$$P_{u}(f)(x)=\sum_{i=0}^{d-1}u_{i}(x)f(\tau_{i}(x)),$$ for all $f :
[0,1] \to \mathbb{R}$ bounded Borel measurable functions.
\end{definition}

A function $h : [0,1] \to \mathbb{R}$ will be called
$P_{u}$-harmonic if $P_{u}(h)=h$ \cite{BJ} \cite{CR}. A  probability  $\nu$ on
$[0,1]$ will be called $P_{u}$-invariant if $P_{u}^{*}(\nu)=\nu$,
where $P_{u}^{*}$ is defined by equality
$$ \int P_{u}(f)(x) d\hat{\nu}=  \int f(x) dP_{u}^{*} \hat{\nu},$$
for all $f : [0,1] \to \mathbb{R}$ continuous.

The correct approach to analyze  an IFS \cite{Elton} \cite{CR}
\cite{BJ} \cite{Ja1} \cite{Ja2} with probabilities $([0,1], \tau_{i}, u_{i})$, is to
consider  for each $x \in [0,1]$, the sequence of random variables
$(Z_{n}(x,.):\Sigma \to [0,1])_{n \in \mathbb{N}}$ as a
realization of the Markov process associated to the Markov chain
with initial distribution $\delta_{x}$ and transitions of
probability $P_{u}$. Moreover, we have a probability
$\mathbb{P}_{x}$ in the space of paths, $\Sigma$, given by
$$\int_{\Sigma}g(w)d\mathbb{P}_{x}= \sum_{i_{1}, ..., i_{n}}
u_{i_{1}}(x)u_{i_{2}}(\tau_{i_{1}}(x))...
u_{i_{n}}(\tau_{i_{n}}...\tau_{i_{1}}(x)) g(i_{1}, ..., i_{n}),$$
when $g=g(x,w)$ depends only of the $n$ first coordinates (see \cite{BJ},
for a proof of the Komolgorov consistence condition).

The probability on path space and the transference operator are
connected by
$$\int_{\Sigma} g(x,w) d\mathbb{P}_{x}(w)=P_{u}^{n}(f)(x),$$
when $g(x,w)=f(\tau_{w_{n}}...\tau_{w_{1}}(x))$ for some
continuous $f$.

\begin{definition}\label{RweightedSystem}
A $\rho$-weighted system, $\rho \geq 0$, is a weighted system
$([0,1], \tau_{i}, u_{i})$ such that there exists a positive
bounded function $h : [0,1] \to
\mathbb{R}$ and $\nu $ probability satisfying\\
$$ P_{u}(h)=\rho h, \quad P_{u}^{*}(\nu)=\rho \nu.$$
\end{definition}

Note that a IFS with probabilities is a $1$-weighted system
(see \cite{DutkayJorg}, \cite{Stenflo} or \cite{BarnDemElton} for
the existence of $P_{u}$-invariant probabilities and
\cite{Stenflo}, Theorem 4, or \cite{StenfloKarlin} for
non-uniqueness of this probabilities). Also, a weighted system,
$([0,1], \tau_{i}, u_{i})$ were $u_{0}=...=u_{d-1}=k$ (constant)
is a $d\, k$-weighted system. Thus the set of $\rho$-weighted
systems is as big class of weighted systems.

Examples of nontrivial $\rho$-weighted system (and
non-probabilistic) can be found in \cite{FanAi-Hua} and
\cite{Stenflo} Corollary 2.

Moreover, from a $\rho$-weighted system $([0,1], \tau_{i}, u_{i})$
one can get
 a normalization $([0,1], \tau_{i}, v_{i})$, in the following way\\
$$ v_{i}(x)= \frac{u_{i}(x) h(\tau_{i}(x))}{\rho h(x)}, \quad \mu=h \nu.$$
Then  $ P_{v}(1)=1$ and $ P_{v}^{*}(\mu)=\mu$.

We thanks an anonymous referee for some comments on a previous version of the present paper. We would like to point out that there exists  some similarities of sections 1, 2, 3 and 5 of our paper with some results in \cite{Ja1} and \cite{Ja2}. We would like to stress that we consider here the holonomic setting which can not be transfer  by some coding to the usual shift case (see remark  4 in section 7). We introduce for such class of probabilities in IFS (which is different from the set of invariant probabilities for $\hat{\sigma}$) the concept of entropy and pressure. It is not the same same concept of entropy as for a measure invariant for the shift $\hat{\sigma}$ (see remarks 3 and 4 in section 7). Also, in our setting, it is natural to consider the all set of possible potentials $u$. In this way our results are of different nature than the ones in \cite{Ja1} \cite{Ja2} where
the dynamical concepts are mainly consider for the shift $\hat{\sigma}$ acting on $  [0,1]\times \Sigma $.

In sections 1 to 6 we consider the basic definitions and results. In sections 7 and 8 we introduce entropy and pressure for holonomic probabilities of IFS.

\vspace{0.3cm}

\section{ Holonomic probabilities}

For IFS we introduce
the concept of \textit{holonomic} probability on $[0,1]\times
\Sigma$ (see \cite{GL} for general definitions and properties in
the setting of symbolic dynamics  of the two-sided shift  in
$\Sigma\times \Sigma$). Several results presented in \cite{GL} can
be easily translated for the IFS setting. In \cite{GL} the main
concern was maximizing probabilities. Here we are mainly
interested in the variational principle of pressure.

By the other hand, some of the new results we presented here can
also be translated to that setting.

\begin{definition}\label{Holonomic}
A holonomic probability  $\hat{\nu}$ on $[0,1]\times \Sigma$ is a
probability such that
$$ \int f(\tau_{X_{1}(w)}(x)) d\hat{\nu}=  \int f(x) d\hat{\nu},$$
for all $f : [0,1] \to \mathbb{R}$ continuous.
\end{definition}

Then the set of holonomic probabilities can be viewed as the set
of probabilities on $[0,1]\times \Sigma$ such that
$$ \int g \circ \hat{\sigma} d\hat{\nu}=  \int g d\hat{\nu},
\quad \forall g \in C([0,1]).$$

From this point of view it is clear that the set of holonomic
probabilities is bigger (see section 4) than the set of $\hat{\sigma}$-invariant
probabilities (because $C([0,1])$ can be viewed as a subset of
$C([0,1]\times \Sigma)$).

\section{Characterization of holonomic probabilities}

Disintegration of probabilities for IFS have been previously consider but for a different purpose \cite{DutJorgDisint}, \cite{DutkayJorg}, \cite{Ja2}.

\begin{definition}\label{RadonSpaces}
A Hausdorff space is a Radon space if all probabilities in this
space is Radon (See~\cite{Schw}).
\end{definition}

\begin{theorem}(\cite{Schw}, Prop. 6, Pg. 117)\label{RadonCaracter}
All compact metric space is Radon.
\end{theorem}

Therefore, all spaces considered here are Radon
spaces.

\begin{theorem}(\cite{Dellach}, Pg 78, (70-III) or \cite{Ambro}, Theorem 5.3.1)\label{Disintegration}
Let $\mathbb{X}$ and $X$ be a separable metric Radon spaces,
$\hat{\mu}$ probability on $\mathbb{X}$, $\pi: \mathbb{X} \to X$
Borel mensurable and $\mu=\pi_{*}\hat{\mu}$. Then there exists a
Borel family of probabilities $\{\hat{\mu}\}_{x \in X}$ on
$\mathbb{X}$, uniquely determined
$\mu$-a.e, such that,\\
1) $\hat{\mu}_{x}(\mathbb{X} \backslash \pi^{-1}(x))=0$, $\mu$-a.e;\\
2) $\int g(z) d\hat{\mu}(z)= \int_{X}\int_{\pi^{-1}(x)}
g(z)d\hat{\mu}_{x}(z) d\mu(x)$.
\end{theorem}
This decomposition is called the \textit{disintegration} of the
probability $\hat{\mu}$.

\begin{theorem}(Holonomic Disintegration)\label{HolonomicDisintegration}
  Consider  a holonomic probability
$\hat{\nu}$ on $[0,1]\times \Sigma$. Let $$\int g(x,w) d\hat{\nu}=
\int_{[0,1]}\int_{\{y\} \times \Sigma} g(x,w) d\hat{\nu}_{y}(x,w)
d\nu(y), \quad \forall g,$$ be the disintegration given by
Theorem~\ref{Disintegration}. Then $\nu$ is  $P_{u}$-invariant for
the IFS with probabilities $([0,1], \tau_{i}, u_{i})_{i=0..d-1}$,
were the  $u_{i}$'s  are given by,
$$u_{i}(y)=\hat{\nu}_{y}(y,\bar{i}), \quad i=0,...,d-1.$$
\end{theorem}
\begin{proof}
Consider a continuous function $f : [0,1] \to \mathbb{R}$ and
defines $I_{1}=\int f(\tau_{X_{1}(w)}(x)) d\hat{\nu}$ and
$I_{2}=\int f(x) d\hat{\nu}$. As  $\hat{\nu}$ is holonomic we have,
$I_{1}=I_{2}$.\\

Now applying the disintegration for both integrals we get\\
$I_{1}=\int_{[0,1]}\int_{\{y\} \times \Sigma}
f(\tau_{X_{1}(w)}(x)) d\hat{\nu}_{y}(x,w) d\nu(y)= \\
=\sum_{i=0}^{d-1}\int_{[0,1]}\int_{\{y\}
\times \bar{i}} f(\tau_{X_{1}(w)}(x)) d\hat{\nu}_{y}(x,w) d\nu(y)=\\
 =\sum_{i=0}^{d-1}\int_{[0,1]} f(\tau_{\bar{i}}(y)) \hat{\nu}_{y}(\{y\}
 \times \bar{i}) d\nu(y)=\int_{[0,1]}P_{u}(f)(y) d\nu(y)$\\
when  $u_{i}(y)=\hat{\nu}_{y}(y,\bar{i}), \quad i=0,...,d-1.$\\

On the other hand\\
$I_{2}=\int_{[0,1]}\int_{\{y\} \times \Sigma} f(x)
d\hat{\nu}_{y}(x,w) d\nu(y)= \int_{[0,1]} f(y) d\nu(y)$.\\

Then, $\int_{[0,1]}P_{u}(f)(y) d\nu(y)= \int_{[0,1]} f(y) d\nu(y)$
for all continuous function $f : [0,1] \to \mathbb{R}$, that is,
$\nu$ is $P_{u}$-invariant.

\end{proof}

\section{Invariance of Holonomic probabilities on IFS}

As we said before, holonomic probabilities are not necessarily
invariant for the map $\hat{\sigma}$. On the other hand all
$\hat{\sigma}$-invariant probability is holonomic. Now we show an
example of holonomic probability which is not $\hat{\sigma}$-invariant (see
\cite{GL} for the case of the two sided  shift).

Suppose that $x_{0} \in [0,1]$, is such that
$Z_{n}(x_{0},\bar{w})=x_{0}$, for some $\bar{w} \in \Sigma$, $n
\in \mathbb{N}$. Then, one can obtain a holonomic probability in
the following way
$$\hat{\nu}=\frac{1}{n} \sum_{j=0}^{n-1} \delta_{\sigma^{j}(\bar{w})}
\times \delta_{Z_{j}(x_{0}, \bar{w})}.$$

Then,
$$\int g(x,w) d\hat{\nu}=\frac{1}{n} \sum_{j=0}^{n-1} g(Z_{j}(x_{0},
\bar{w}) ,\sigma^{j}(\bar{w})),$$

Note that this probability is holonomic but not
$\hat{\sigma}$-invariant.

In fact, it is enough to see that
$$\int g \circ \hat{\sigma}(x,w)
d\hat{\nu}=\frac{1}{n} \sum_{j=0}^{n-1} g \circ
\hat{\sigma}(Z_{j}(x_{0}, \bar{w}) ,\sigma^{j}(\bar{w}))=$$
$$=\frac{1}{n} \sum_{j=0}^{n-1} g (Z_{j+1}(x_{0}, \bar{w})
,\sigma^{j+1}(\bar{w})).$$

Thus, $\int g \circ \hat{\sigma}(x,w) d\hat{\nu}- \int g(x,w)
d\hat{\nu}=\frac{1}{n} g(x_{0}, \sigma^{n}(\bar{w}))- g(x_{0},
\bar{w}))$, and it is clearly not identical to 0, $\forall g$.

However, $\hat{\nu}$ is holonomic because given any continuous
function $f : [0,1] \to \mathbb{R}$ we have
$$\int f(\tau_{X_{1}(w)}(x)) d\hat{\nu}=\frac{1}{n} \sum_{j=0}^{n-1}
f(\tau_{X_{1}(\sigma^{j}(\bar{w}) )}(Z_{j}(x_{0}, \bar{w}))=$$
$$=\frac{1}{n} \sum_{j=0}^{n-1} f(Z_{j+1}(x_{0}, \bar{w}))=\int f(x) d\hat{\nu},$$
because,  $Z_{n}(x_{0},\bar{w})=x_{0}$.
\vspace{0.3cm}

\section{Ergodicity of holonomic probabilities}

Given a holonomic probability $\hat{\nu}$, we can associate, by
holonomic disintegration, a unique IFS with probabilities $([0,1],
\tau_{i}, u_{i})$ such that $P_{u}^{*}(\nu)=\nu$ and
$\nu=\pi_{*}\hat{\nu}$.

Let $Z_{n}(\cdot) :[0,1] \hookleftarrow, \; n\in  \mathbb{N}$, be a
sequence of random variables on $[0,1]$. Then, we obtain a
Markov process with transition of probabilities $P_{u}$ and initial
distribution $\nu$, that we will denote by $(Z_{n}, P_{u}, \nu)$.

This process is a stationary process by construction, thus does make
sense to ask if $(Z_{n}, P_{u}, \nu)$ is ergodic(\cite{Elton} for
details of this process and definition of ergodicity).

\begin{definition}\label{HolonomicErgodic}
A holonomic probability $\hat{\nu}$ is called ergodic, if the
associated Markov process $(Z_{n}, P_{u}, \nu)$ is an ergodic
process.
\end{definition}

\begin{lemma}\label{HolonUniqErgodic}(Elton,\cite{Elton})
   Let $\hat{\nu}$ be a holonomic probability with holonomic
   disintegration $([0,1], \tau_{i}, u_{i})$. If $\pi_{*}\hat{\nu}$
   is the unique $P_{u}$-invariant probability, then  $\hat{\nu}$ is
   ergodic.
\end{lemma}

\begin{theorem}\label{HolonErgodicTheorem}(Elton,\cite{Elton})
   Let $([0,1], \tau_{i})$ be a contractive IFS (contractiveness means
   that $\tau_{i}$ is a contraction for all $i$) and $\hat{\nu}$ be
   a ergodic holonomic probability with holonomic
   disintegration $([0,1], \tau_{i},
   u_{i})$. Suppose that $u_{i} \geq \delta >0, \; \forall
   i=0,...,d-1$. Then, $\forall x \in [0,1]$ there exists $G_{x}
   \subseteq \Sigma$ such that $\mathbb{P}_{x}(G_{x})=1$ and for
   each $w \in G_{x}$
   $$      \frac{1}{N}\sum_{i=0}^{N-1}       f(Z_n (x,w))  =     \frac{1}{N}\sum_{i=0}^{N-1}f(\hat{\sigma}^{i}(x,w)) \to \int f
   d\hat{\nu}= \int f(x) d \hat{\nu}(x,w),$$
for all $f \in C([0,1])$.
\end{theorem}
\begin{proof}
The proof is a straightforward modification of the one presented in Elton's ergodic
theorem (see \cite{Elton} \cite{Ja2}). In fact, the contractiveness of
$([0,1], \tau_{i})$ its stronger that Dini condition that appear in
Elton's proof (see \cite{Elton} and \cite{BarnDemElton}) and the
ergodicity of $\hat{\nu}$ can replace the  uniqueness  of the
initial distribution in the last part of the argument. The other
parts of the proof are the same as in \cite{Elton}.
\end{proof}

We point out that Elton's Theorem is not the classical ergodic theorem for $\hat{\sigma}$. The claim is:
$\forall x \in [0,1]$ there exists $G_{x}\subseteq \Sigma$ such that $\mathbb{P}_{x}(G_{x})=1$ and for each $w \in G_{x}$ ... Moreover, $f:[0,1] \to \mathbb{R}$.

   This theorems fits well for holonomic probabilities in the IFS case. We just mention it in order to give to the reader a broader perspective of the holonomic setting. We do not use it in the rest of the paper.

\section{Construction of holonomic probabilities
for $\rho$-weighted systems}

Given a  $\rho$-weighted system $([0,1], \tau_{i}, u_{i})$, that
is, $$ P_{u}(h)=\rho h, \quad P_{u}^{*}(\nu)=\rho \nu,$$ consider
the normalization  $([0,1], \tau_{i}, v_{i})$, then $ P_{v}(1)=1$
and $ P_{v}^{*}(\mu)=\mu$.

Its easy to see that the probability on $[0,1]\times \Sigma$ given
by
$$\int g(x,w) d\hat{\mu}= \int_{[0,1]}\int_{\Sigma} g(x,w)d\mathbb{P}_{x}(w)
d\mu(x)$$ is holonomic if $\mathbb{P}_{x}$ is given from $v$ (see
\cite{DutJorgDisint}\cite{Ja1}\cite{Ja2} for disintegration of projective measures on
IFS). The probability $\hat{\mu}$ will be called the
\textit{holonomic lifting} of $\mu$.

{\bf Remark 1}. We point out that  that the holonomic lifting $\hat{\mu}$ of a given $\mu$ (as above) is a $\hat{\sigma}$-invariant probability (one can see that by taking functions that depends only of finite symbols and applying de definition of a ${P}_{x}$ probability). So $$\pi_{*} \{\emph{Holonomic} \; probabilities \}= \pi_{*} \{ \hat{\sigma}-invariant \; probabilities \}.$$

We will consider in the next sections the concept of pressure. The value of pressure among holonomic or among invariant will be the same.
However one cannot reduce the study of variational problems involving holonomic probabilities to the study of $\hat{\sigma}$-invariant probabilities. This will be explained in remark 4 in example 3 after Theorem~\ref{VarPrinciple}.

One can reverse the process, starting from a IFS with
probabilities (a 1-weighted system) $([0,1], \tau_{i}, v_{i})$,
that is, $ P_{v}(1)=1$ and $ P_{v}^{*}(\nu)=\nu$, and consider the
associated $\hat{\nu}$, the holonomic lifting of $\nu$.
$$\int g(x,w) d\hat{\nu}= \int_{[0,1]}\int_{\Sigma} g(y,w)d\mathbb{P}_{y}(w)
d\nu(y).$$

By holonomic
disintegration (Theorem~\ref{HolonomicDisintegration}), one can
represents the probability $\hat{\nu}$ as
$$\int g(x,w) d\hat{\nu}=
\int_{[0,1]}\int_{\{y\} \times \Sigma} g(x,w) d\hat{\nu}_{y}(x,w)
d\nu_{0}(y).$$

Then, $\nu_{0}$ is  $P_{u}$-invariant for the IFS
with probabilities $([0,1], \tau_{i}, u_{i})_{i=0..d-1}$, were the
$u_{i}$'s  are given by,
$$u_{i}(y)=\hat{\nu}_{y}(y,\bar{i}), \quad
i=0,...,d-1\,.$$

We point out  that $\nu_{0}= \nu$ (it is a straightforward
calculation), moreover we can rewrite
$$\int g(x,w) d\hat{\nu}= \int_{[0,1]}\int_{y \times \Sigma} g(x,w)
d(\delta_{y}(x) \times \mathbb{P}_{y}(w)) d\nu(y).$$

By the uniqueness in Theorem~\ref{Disintegration} we get,
$$\hat{\nu}_{y}=\delta_{y} \times \mathbb{P}_{y}, \quad \nu-a.e\,.$$

Then, we have
$$u_{i}(y)=\hat{\nu}_{y}(y,\bar{i})=(\delta_{y} \times
\mathbb{P}_{y})(y,\bar{i})=\mathbb{P}_{y}(\bar{i})=v_{i}(x),
\quad \nu-a.e\,.$$

From this argument we get the following proposition
 \begin{proposition} \label{Lifting}
 Let $([0,1], \tau_{i}, v_{i})$ be a 1-weighted system and
 $\hat{\nu}$ the holonomic lifting of the invariant probability $P_{v}$-invariant $\nu$.
 If $([0,1], \tau_{i}, u_{i})$ is the  1-weighted system obtained
 from holonomic disintegration of $\hat{\nu}$, then $u=v, \quad
 \nu-a.e$, where $\nu=\pi_{*}\hat{\nu}$.
 \end{proposition}

\vspace{0.3cm}

\section{Entropy and a variational principle for $\rho$-weighted systems}

   Let us consider a $\rho$-weighted system, $([0,1], \tau_{i}, u_{i})$.
   Denote by $\mathbb{B}^{+}$ the set of all positive bounded Borel functions
   on $[0,1]$ and by $\mathcal{H}$ the set of all holonomic probabilities on $[0,1]\times \Sigma$ for
   $\hat{\sigma}$.

The central idea in this section is to consider a generalization
of the definition of entropy for the case of holonomic
probabilities via the concept naturally suggested by Theorem 4 in
\cite{ALopes}. We will show that under such point of view the
classical results in Thermodynamic Formalism are also true.

   Given $\hat{\nu} \in \mathcal{H}$ we can define the functional
   $\alpha_{\hat{\nu}}: \mathbb{B}^{+} \to \mathbb{R} \cup \{-\infty\}$
   by  $$ \alpha_{\hat{\nu}}(\psi)=\inf_{f \in \mathbb{B}^{+}} \int
 \ln(\frac{P_{\psi}f}{\psi f}) d\hat{\nu}.$$

    Let $\alpha_{\hat{\nu}}$ be the functional defined above.
    Observe that  $\alpha_{\hat{\nu}}$ doesn't depend of $\psi$.

    In fact, by taking $\psi_{1}, \psi_{2} \in \mathbb{B}^{+}$ and
    $f \in \mathbb{B}^{+}$, define $ g  \in \mathbb{B}^{+}$ as
    being $g=\frac{\psi_{1}}{\psi_{2}} f \in  \mathbb{B}^{+}$,
    then $$\int \ln(\frac{P_{\psi_{2}}g}{\psi_{2} g}) d\hat{\nu}=\int
 \ln(\frac{P_{\psi_{1}}f}{\psi_{1} f}) d\hat{\nu} .$$
    Thus, $\alpha_{\hat{\nu}}(\psi_{2})=\alpha_{\hat{\nu}}(\psi_{1})$.

   \begin{definition}\label{Entropy}
    Given $\hat{\nu} \in \mathcal{H}$ we  define the Entropy
    of $\hat{\nu}$ by $$h(\hat{\nu})=\alpha_{\hat{\nu}}.$$
   \end{definition}

   From above we get
   $$h(\hat{\nu})=\inf_{f \in \mathbb{B}^{+}} \int \ln(\frac{P_{\psi}f}{\psi f}) d\hat{\nu},$$
   $\forall \psi \in \mathbb{B}^{+}$.

The above definition agrees with the usual one for invariant probabilites when it is consider a transformation of degree $d$, its $d$-branches and the naturally associated IFS (see \cite{ALopes}).

      \begin{lemma} \label{LogaritInequal}
   Given $\beta \geq  1 + \alpha$ and numbers $a_{i} \in [1+ \alpha ,\beta], \;
   i=0,...,d-1$ there exists $\varepsilon \geq 1$ such that
   $$ \ln(\varepsilon \sum_{i=0}^{d-1} a_{i}) =  \sum_{i=0}^{d-1} \ln(\varepsilon
   a_{i}).$$
   \end{lemma}

   This lemma follows from the choice $\varepsilon =
   exp(\frac{1}{d-1} ( \ln(\sum_{i=0}^{d-1} a_{i}))/ (\sum_{i=0}^{d-1}
 \ln(a_{i}))$ and the fact that $a_{i} \geq 1+ \alpha $.

   \begin{lemma} \label{HolonomicInequal}
   Given $f \in \mathbb{B}^{+}$ and  $\hat{\nu} \in \mathcal{H}$ then
   $$ \sum_{i=0}^{d-1} \int  f(\tau_{i}(x)) d\hat{\nu} \geq   \int  f(x) d\hat{\nu}.$$
.   \end{lemma}
  \begin{proof}
   As  $\hat{\nu}$ is holonomic, then we have
   $$\int f(\tau_{X_{1}(w)}(x)) d\hat{\nu}=
   \int f(x) d\hat{\nu}.$$
   This can be written  as
   $$\int f(\tau_{X_{1}(w)}(x)) d\hat{\nu}=\sum_{i=0}^{d-1}
   \int_{[0,1] \times \overline{i}} f(\tau_{X_{1}(w)}(x)) d\hat{\nu} =\sum_{i=0}^{d-1}
   \int_{[0,1] \times \overline{i}} f(\tau_{i}(x)) d\hat{\nu} .$$

   Note that for each $i$
   $$\int f(\tau_{i}(x)) d\hat{\nu}=\sum_{j=0}^{d-1}
   \int_{[0,1] \times \overline{j}} f(\tau_{i}(x)) d\hat{\nu} \geq
   \int_{[0,1] \times \overline{i}} f(\tau_{i}(x)) d\hat{\nu} .$$

   Thus,
   $$\sum_{i=0}^{d-1}
   \int f(\tau_{i}(x)) d\hat{\nu} \geq \sum_{i=0}^{d-1}
   \int_{[0,1] \times \overline{i}} f(\tau_{i}(x)) d\hat{\nu} = \int  f(x) d\hat{\nu}.$$
   \end{proof}

   \begin{proposition} \label{PositOfEntropy}
     Consider $\hat{\nu} \in \mathcal{H}$. Then $$0 \leq h(\hat{\nu}) \leq \ln(d).$$
   \end{proposition}
   \begin{proof}
   Initially, consider $\psi=1$. We know that
   $h(\hat{\nu})=\inf_{f \in \mathbb{B}^{+}} \int \ln(\frac{P_{1}f}{f}) d\hat{\nu} \leq \int \ln(\frac{P_{1}1}{1}) d\hat{\nu}= \ln(d)$.

   Now, in order to prove the inequality  $$h(\hat{\nu}) =\inf_{f \in \mathbb{B}^{+}}
   \int \ln(\sum_{i=0}^{d-1} \frac{ f \circ \tau_{i}}{f}) d\hat{\nu} \geq 0,$$
   consider $I=\int \ln(\sum_{i=0}^{d-1} \frac{ f \circ \tau_{i}}{f})
   d\hat{\nu}$ and suppose, without lost of generality,  that $1 + \alpha \leq f
   \leq \beta$ (because this integral is invariant under the projective function $f \to \lambda f$).
   Then, we can write this integral as
   $$I= \int \ln(\sum_{i=0}^{d-1} \frac{\varepsilon f \circ \tau_{i}}{\varepsilon f})
   d\hat{\nu}= \int \ln(\varepsilon \sum_{i=0}^{d-1}  f \circ \tau_{i})
   d\hat{\nu} -\int \ln(\varepsilon f) d\hat{\nu}, \quad \quad  (1)$$
   In order to use the inequality obtained from Lemma~\ref{LogaritInequal},
   denote (for each fixed $x$) $$a_{i}= f \circ \tau_{i}(x).$$
   Then, we get
   $$\varepsilon(x) =
   exp(\frac{1}{d-1} \cdot \frac{ \ln(\sum_{i=0}^{d-1}f \circ \tau_{i})}{ (\sum_{i=0}^{d-1}
 \ln(f \circ \tau_{i}))}) \geq \varepsilon_{0} \geq 1,$$ by the
   compactness of $[0,1]$.
   From this choice we get
   $$ \ln(\varepsilon_{0} \sum_{i=0}^{d-1} f \circ \tau_{i}) \geq
   \sum_{i=0}^{d-1} \ln(\varepsilon_{0} f \circ \tau_{i}). \quad \quad (2)$$
   Using (2) in (1) we obtain
   $$I \geq \sum_{i=0}^{d-1} \int \ln(\varepsilon_{0}  f \circ \tau_{i})
   d\hat{\nu} -\int \ln(\varepsilon_{0} f) d\hat{\nu}.$$
   Moreover,  using the inequality from Lemma~\ref{HolonomicInequal} applied to the
   function $ \ln(\varepsilon f)$ (note that $\ln (\varepsilon_{0} f) \in
   \mathbb{B}^{+}$, because $\varepsilon_{0} \geq 1$), we get
   $$I \geq \int \ln(\varepsilon f) d\hat{\nu} -\int \ln(\varepsilon f) d\hat{\nu} =0.$$
   \end{proof}

   \begin{definition}\label{Pressure}
    Given $\phi \in \mathbb{B}^{+}$ we  define the Topological Pressure
    of $\phi$ by $$p(\phi)= \sup_{\hat{\nu} \in \mathcal{H}} \{ h(\hat{\nu}) +
    \int \ln(\phi) d\hat{\nu} \}.$$
   \end{definition}

{\bf Remark 2}. From Remark 1 it follows that
$$\sup_{\hat{\nu} \in \mathcal{H}} \{ h(\hat{\nu}) +
    \int \ln(\phi) d\hat{\nu} \}=$$
    $$ \sup_{\hat{\nu} \,\text{invariant for} \, \hat{\sigma}} \{ h(\hat{\nu}) +
    \int \ln(\phi) d\hat{\nu} \}.$$

\vspace{0.2cm}

   Using the formula for entropy we get a characterization of topological pressure as
   $$p(\phi)= \sup_{\hat{\nu} \in \mathcal{H}} \{ \inf_{f \in \mathbb{B}^{+}}
    \int \ln(\frac{P_{\phi}f}{\phi f}) d\hat{\nu} +  \int \ln(\phi) d\hat{\nu} \} =
     \sup_{\hat{\nu} \in \mathcal{H}} \inf_{f \in \mathbb{B}^{+}}
      \int \ln(\frac{P_{\phi}f}{f}) d\hat{\nu}.$$

   \begin{definition}\label{EquilibriumState}
    A holonomic measure $\hat{\nu}_{eq}$ such that
    $$p(\phi)= h(\hat{\nu}_{eq}) + \int \ln(\phi) d\hat{\nu}_{eq}$$
    will be called an equilibrium state for $\phi$.
   \end{definition}

   {\bf Remark 3}. Example 3 bellow shows that in IFS there exist examples of holonomic equilibrium states for $\phi$ which are not invariant for $\hat{\sigma}.$
   \vspace{0.2cm}

   In the next theorem we do not assume
   $\sum_{i=0}^{d-1}\phi(\tau_{i}(x))=1$.

\begin{theorem} \label{PressureIsLogRo}
  Let us consider $\phi \in \mathbb{B}^{+}$ such that $([0,1], \tau_{i}, \phi)$ is
  a $\rho$-weighted system, for some $\rho \geq 0$. Then, $p(\phi)= \ln(\rho)$.
  In particular,  the transference operator $P_{\phi}$ has a unique positive eigenvalue.
\end{theorem}
\begin{proof}
     Note  that, $h(\hat{\nu})=\inf_{f \in \mathbb{B}^{+}}
     \int \ln(\frac{P_{\phi}f}{\phi f}) d\hat{\nu} \leq \int \ln(\frac{P_{\phi}h}{\phi h})
     d\hat{\nu}= \int \ln(\frac{\rho}{\phi}) d\hat{\nu}= -\int \ln(\phi) d\hat{\nu} + \ln(\rho)$
     so, $$h(\hat{\nu}) +\int \ln(\phi) d\hat{\nu} \leq \ln(\rho), \; \forall \hat{\nu},$$
     thus, $p(\phi) \leq \ln(\rho)$.

     Remember that $$p(\phi)=\sup_{\hat{\nu} \in \mathcal{H}}
     \inf_{f \in \mathbb{B}^{+}} \int \ln(\frac{P_{\phi}f}{f}) d\hat{\nu}.$$
     Let $\hat{\nu}_{0}$ be a fixed  holonomic probability  such that the normalized
      dual operator verifies $P_{u}^{*} (\pi_{*}\hat{\nu}_{0})=\pi_{*}\hat{\nu}_{0}$
      (always there exists if $P_{u}$ is the normalization of $P_{\phi}$),
      were $\pi(x,\omega)=x$.  Thus we can write
     $$p(\phi) \geq  \inf_{f \in \mathbb{B}^{+}} \int \ln(\frac{P_{\phi}f}{f})
     d\hat{\nu}_{0}.$$ Note that, from the normalization property we get
     $$P_{\phi}(f)=P_{u}(\frac{f}{h}) \rho h, \quad \forall f.$$

     Moreover, we know that $\ln(P_{u}g) \geq P_{u}\ln(g), \forall g$, by concavity
     of logarithmic function.

     Now, considering an arbitrary $f \in \mathbb{B}^{+}$, we get
     $$\int \ln(\frac{P_{\phi}f}{f}) d\hat{\nu}_{0}=
     \int \ln(\frac{P_{u}(f/h)\rho h}{f}) d\hat{\nu}_{0}=\int \ln(\frac{P_{u}(f/h)}{f/h})
     d\hat{\nu}_{0} + \ln(\rho) \geq $$
      $$ \geq \int P_{u} \ln(f/h) d\hat{\nu}_{0}- \int \ln(f/h) d\hat{\nu}_{0} + \ln(\rho) =
 \ln(\rho).$$

 So, $\inf_{f \in \mathbb{B}^{+}} \int \ln(\frac{P_{\phi}f}{f})
       d\hat{\nu}_{0} \geq \ln(\rho)$, that is, $p(\phi) \geq \ln(\rho)$. From
       this we get  $p(\phi)= \ln(\rho)$.

     In order to obtain the second part of the claim it is enough to see that $p(\phi)= \ln(\rho)$,
     for all $\rho$, thus the eigenvalue is unique.
\end{proof}

\begin{theorem}(Variational principle)\label{VarPrinciple}
  Consider $\phi \in \mathbb{B}^{+}$ such that $([0,1], \tau_{i}, \phi)$ is a
  $\rho$-weighted system, for $\rho = e^{p(\phi)} \geq 0$. Then, any holonomic probability
   $\hat{\nu}_{0}$ such that the normalized dual operator verifies
    $P_{u}^{*} (\pi_{*}\hat{\nu}_{0})=\pi_{*}\hat{\nu}_{0}$ is an equilibrium state.
\end{theorem}
\begin{proof}
    Note that, from the definition of pressure, we get
    $\ln(\rho)=p(\phi) \geq h(\hat{\nu}_{0}) + \int \ln(\phi) d\hat{\nu_{0}}$.
     As, $P_{u}^{*} (\pi_{*}\hat{\nu}_{0})=\pi_{*}\hat{\nu}_{0}$,
     for an arbitrary $f \in \mathbb{B}^{+}$, we obtain
    $$h(\hat{\nu}_{0}) + \int \ln(\phi) d\hat{\nu_{0}}=\int \ln(\frac{P_{\phi}f}{f})
     d\hat{\nu}_{0} \geq \ln(\rho).$$
    Thus, $\ln(\rho)= h(\hat{\nu}_{0}) + \int \ln(\phi) d\hat{\nu_{0}}$.

\end{proof}

    Note that, if $\hat{\nu}_{0}$ is a equilibrium state, such that
    the IFS with probabilities $([0,1], \tau_{i}, v_{i})$ associated
    to $\hat{\nu}_{0}$ by holonomic disintegration, is uniform, then
    $P_{u}^{*} (\pi_{*}\hat{\nu}_{0})=\pi_{*}\hat{\nu}_{0}$.

    In fact, we know that $\ln(\rho)= h(\hat{\nu}_{0}) + \int \ln(\phi)
    d\hat{\nu_{0}}$, and $P_{v}^{*}
    (\pi_{*}\hat{\nu}_{0})=\pi_{*}\hat{\nu}_{0}$. Then, we can write
    $$\ln(\rho)= h(\hat{\nu}_{0}) + \int P_{v} \ln(\phi)
    d\hat{\nu_{0}}= h(\hat{\nu}_{0}) + \int \sum_{i=0}^{d-1}
    v_{i} \ln( \phi(\tau_{i})) d\hat{\nu}_{0} .$$

    Remember that the normalization of $\phi$ is given by
    $$ u_{i}(x)= \frac{\phi(\tau_{i}(x)) h(\tau_{i}(x))}{\rho
    h(x)},$$ thus,
    $$0= h(\hat{\nu}_{0}) + \int \sum_{i=0}^{d-1}
    v_{i} \ln( u_{i}) d\hat{\nu}_{0}. \quad \quad (3)$$

    As $([0,1], \tau_{i}, v_{i})$ is uniform, there exists a
    weight function $\psi$ such that $v_{i}(x)=\phi(\tau_{i}(x)), \; \forall
    i=0,...,d-1$. Moreover, $p(\psi)=0$ and $\hat{\nu_{0}}$ is
    clearly a equilibrium state for $\psi$. Thus
    $$0= h(\hat{\nu}_{0}) + \int \ln(\psi)
    d\hat{\nu_{0}}$$
    Using $P_{v}^{*} (\pi_{*}\hat{\nu}_{0})=\pi_{*}\hat{\nu}_{0}$ we
    get
    $$0= h(\hat{\nu}_{0}) + \int \sum_{i=0}^{d-1}
    v_{i} \ln( v_{i}) d\hat{\nu}_{0} .\quad \quad (4)$$

    It is well known that\\
    $$-\sum_{i=0}^{d-1}  a_{i} \ln( a_{i}) + \sum_{i=0}^{d-1}
    a_{i} \ln( b_{i}) \leq 0.   \quad \quad (5)$$
    where $\sum_{i=0}^{d-1}  a_{i}=1=\sum_{i=0}^{d-1}  b_{i}$
    and $b_{i} \geq 0$, with equality only if $a_{i}=b_{i}$ (see \cite{Parry} for a
    proof). From  (3) and (4) we get,
    $$u_{i}(x)=v_{i}(x), \quad \nu-a.e.$$
    Then, it follows that $P_{u}^{*}
    (\pi_{*}\hat{\nu}_{0})=\pi_{*}\hat{\nu}_{0}$.
\vspace{0.3cm}

{\bf Examples:}

   1) For $\phi =1$, we have $P_{\phi} (1)= d  \cdot 1 $. Then, for all equilibrium
    states $\hat{\nu}_{eq}$ we get  $$\ln(d) = h(\hat{\nu}_{eq}) +
     \int \ln(1) d\hat{\nu}_{eq}=
   \sup_{\hat{\nu} \in \mathcal{H}} \inf_{f \in \mathbb{B}^{+}}
   \int \ln(\frac{P_{1}f}{f}) d\hat{\nu}= \sup_{\hat{\nu} \in \mathcal{H}} h(\hat{\nu})$$

   2) If $P_{\phi} (1)=  1 $, that is, the case of IFS with probabilities, then $p(\phi)=0$.
    Therefore, $h(\hat{\nu}) + \int \ln(\phi) d\hat{\nu} \leq 0$, for all holonomic probabilities.
    Moreover, any equilibrium state $\hat{\nu}_{eq}$ satisfies
   $$h(\hat{\nu}_{eq})\, = - \, \int \ln(\phi) d\hat{\nu}_{eq}$$

   3) Consider the IFS given by $([0,1], \tau_{0}(x)=x, \tau_{1}(x)=1-x)$ and the potential $\phi(x)=2+\cos(2 \pi x)$. Is clear that $$0 \leq \int \ln(\phi) d\nu \leq \ln 3, \; \forall \nu.$$  We will consider a especial holonomic probability $\hat{\nu}_{0}$ constructed in the following way (similar to the one  presented in section 4):\\

   Consider fixed $\bar{w} \neq (111111...) \in \Sigma$ and $x_{0}=0$. The holonomic probability $\hat{\nu}_{0}$ is the  average of delta of Dirac distributions at the points $(x_{0},11\bar{w}) $ and  $\hat{\sigma} (x_{0},11\bar{w})$  in $[0,1] \times \Sigma$, more precisely, for any $g$
   $$\int g(x,w) d\hat{\nu}=\frac{1}{2} \sum_{j=0}^{1} g(Z_{j}(x_{0},
   \bar{w}) ,\sigma^{j}(\bar{w}))=\frac{1}{2} ( g(0,11\bar{w}) + g(1,1\bar{w}))$$

   This probability is not $\hat{\sigma}$-invariant by construction, and has the interesting property:
   $$\ln(2) = h(\hat{\nu}_{0})= \sup_{\hat{\nu} \in \mathcal{H}} h(\hat{\nu})$$

   Indeed, since $h(\hat{\nu}) \leq \ln(2), \forall \hat{\nu}$, it is enough to see that $h(\hat{\nu}_{0})= \ln(2)$. Remember that $h(\hat{\nu}_{0})= \inf_{f \in \mathbb{B}^{+}} \int \ln(\frac{P_{1}f}{f}) d\hat{\nu}_{0}$, so
   $$h(\hat{\nu}_{0})= \inf_{f \in \mathbb{B}^{+}} \frac{1}{2} ( \log \frac{P_{1}f(0)}{f(0)} + \log \frac{P_{1}f(1)}{f(1)})=\inf_{f \in \mathbb{B}^{+}} \frac{1}{2}[ \ln(1 +\frac{f(0)}{f(1)}) + \ln(1 +\frac{f(1)}{f(0)})]$$

   Taking $\frac{f(0)}{f(1)}= \lambda > 0 $, we get
   $$h(\hat{\nu}_{0})= \inf_{\lambda > 0} \frac{1}{2} [ \ln(1 + \lambda) + \ln( 1 +\frac{1}{\lambda}) ]= \ln 2 $$

  {\bf Remark 4}.  This shows that such $\hat{\nu}_{0}$ is a maximal entropy holonomic probability which is not $\hat{\sigma}$-invariant. This also shows that the holonomic setting can not be reduced, via coding, to the analysis of $\hat{\sigma}$-invariant probabilities in a symbolic space. Otherwise, in the symbolic case a
  probability with support in two points would have positive entropy.

\vspace{0.2cm}

   Also, one can calculate $m=\sup_{\hat{\nu} \in \mathcal{H}} \int \ln(\phi) d\hat{\nu} \leq \ln 3$.
    We claim that $m=\ln 3$. In fact,
   $$\int \ln(\phi) d\hat{\nu}_{0} = \frac{1}{2}[ \ln(2+\cos(2 \pi 0)) + \ln(2+\cos(2 \pi 1))]=\ln 3.$$

   Finally, we point out that $\hat{\nu}_{0}$ is also a equilibrium state. Indeed,

   $$p(\phi)=\sup_{\hat{\nu} \in \mathcal{H}} \{ h(\hat{\nu}) +
    \int \ln(\phi) d\hat{\nu} \} \leq \ln 2 + \ln 3 .$$
   As, $h(\hat{\nu}_{0}) + \int \ln(\phi) d\hat{\nu}_{0}= \ln 2 + \ln 3$, we get
   $$\ln(6)= p(\phi)= h(\hat{\nu}_{0}) + \int \ln(\phi) d\hat{\nu}_{0}.$$

   From this example one can see that there exists equilibrium states which are not $\hat{\sigma}$-invariant probabilities.

   \begin{definition}\label{Cohomology}
    Two functions $\psi_{1},\psi_{2} \in \mathbb{B}^{+}$ will be called
    holonomic-equivalent (or, co-homologous) if there exists a function $h \in \mathbb{B}^{+}$, such that
    $$ \psi_{1}(x) = \psi_{2}(x) \cdot \frac{h(\tau_{X_{1}(w)}(x))}{h(x)}, \; \forall (x,\omega).$$
   \end{definition}

   It is clear that, two holonomic equivalent potentials $\psi_{1},\psi_{2}$ will have the same
   equilibrium states.

\vspace{0.3cm}

\section{An alternative point of view for  the concept of  entropy and pressure for IFS}
   \begin{definition}\label{EntropyAlt}
    Given $\hat{\nu} \in \mathcal{H}$ and let $([0,1], \tau_{i}, v_{i})$  be the IFS
     with probabilities arising in the holonomic disintegration of $\hat{\nu}$ (see
     Theorem~\ref{HolonomicDisintegration}). We can also define the Entropy   of $\hat{\nu}$ by
    $$h(\hat{\nu})= -\sup_{\sum_{i=0}^{d-1} u_{i}=1} \int \sum_{i=0}^{d-1}
    v_{i} \ln( u_{i}) d\hat{\nu}.$$
   \end{definition}

   \begin{proposition} \label{PositOfEntropyAlt}
     Consider $\hat{\nu} \in \mathcal{H}$. Then, $$0 \leq h(\hat{\nu}) \leq \ln(d).$$
   \end{proposition}
   \begin{proof}
   Firstly consider $u_{i}^{0}=\frac{1}{d}, \; i=0,...,d-1$ then
   $\sum_{i=0}^{d-1} u_{i}^{0}=1$ and
   $$-h(\hat{\nu})= \sup_{\sum_{i=0}^{d-1} u_{i}=1} \int \sum_{i=0}^{d-1}
    v_{i} \ln( u_{i}) d\hat{\nu} \geq \int \sum_{i=0}^{d-1}
    v_{i} \ln( \frac{1}{d}) d\hat{\nu}= -\ln(d),$$
   so, $h(\hat{\nu}) \leq \ln(d)$.

    On the other hand,  $u_{i} \leq 1$ so $\ln( u_{i}) \leq 0$, and
    then $$\sup_{\sum_{i=0}^{d-1} u_{i}=1} \int \sum_{i=0}^{d-1}
    v_{i} \ln( u_{i}) d\hat{\nu} \leq 0.$$ Thus, $0 \leq h(\hat{\nu})$.
   \end{proof}

   \begin{lemma} \label{ExistEquilStateAlt}(Existence of equilibrium
   states)   Consider $\phi \in \mathbb{B}^{+}$ such that $([0,1], \tau_{i}, \phi)$ is
  a $\rho$-weighted system, for some $\rho > 0$ (remember that there exists $h > 0$,
  such that $ P_{\phi}(h)=\rho h$). Denote $P_{v}$  the normalization of
  $P_{\phi}$, that is,  $([0,1], \tau_{i}, v_{i})$ is a 1-weighted system,
  such that $ P_{v}(1)=1$  and $ P_{v}^{*}(\nu)=\nu$.
  Let $\hat{\nu}$ be the holonomic lifting of $\nu$. Then,
$$h(\hat{\nu}) + \int \ln(\phi) d\hat{\nu} = \ln(\rho)$$
\end{lemma}
   \begin{proof}
    Let $\hat{\nu}$ be the holonomic lifting of $\nu$. By
    Proposition~\ref{Lifting}, we know that the 1-weighted system
    associated to its holonomic disintegration is $([0,1], \tau_{i},
    v_{i})$. Then, from the definition of entropy we get
    $$h(\hat{\nu})= -\sup_{\sum_{i=0}^{d-1} u_{i}=1} \int \sum_{i=0}^{d-1}
    v_{i} \ln( u_{i}) d\hat{\nu}= - \int \sum_{i=0}^{d-1}
    v_{i} \ln( v_{i}) d\hat{\nu} $$
    from the logarithmic inequality (5) above.

    Remember that the normalization of $\phi$ is given by
    $$ v_{i}(x)= \frac{\phi(\tau_{i}(x)) h(\tau_{i}(x))}{\rho h(x)}.$$

    Replacing this expression in the equation for entropy we get

    $$h(\hat{\nu})= - \int \sum_{i=0}^{d-1} v_{i}
    \ln(\frac{\phi(\tau_{i}(x)) h(\tau_{i}(x))}{\rho h(x)}) d\hat{\nu} =$$
    $$= - \int \sum_{i=0}^{d-1} v_{i}  \ln(\phi(\tau_{i}(x))
    d\hat{\nu}  - \int \sum_{i=0}^{d-1} v_{i}  \ln(\frac{h(\tau_{i}(x))}{h(x)})
    d\hat{\nu} + \ln(\rho)=$$
    $$= - \int \ln(\phi(x))d\hat{\nu} + \ln(\rho)$$
\end{proof}

Now, we use the concept introduced in the present section.

\begin{definition}\label{PressureAlt}
    Given $\phi \in \mathbb{B}^{+}$, we  define the Topological Pressure
    of $\phi$ by $$p(\phi)= \sup_{\hat{\nu} \in \mathcal{H}} \{ h(\hat{\nu}) +
     \int \ln(\phi) d\hat{\nu} \}$$
\end{definition}

\begin{theorem} \label{PressureIsLogRoAlt}
  Let us consider $\phi \in \mathbb{B}^{+}$ such that $([0,1], \tau_{i}, \phi)$ is
  a $\rho$-weighted system, for some $\rho \geq 0$. Then $p(\phi)=\ln(\rho)$.
  In particular,  the transference operator $P_{\phi}$ has a unique positive eigenvalue.
\end{theorem}

\begin{proof}
Let $P_{u}$ be the normalization of $P_{\phi}$. Then,
$$h(\hat{\nu}) + \int \ln(\phi) d\hat{\nu} = h(\hat{\nu}) + \int \ln(\phi) d\hat{\nu} =$$
$$= - \int \sum_{i=0}^{d-1}  v_{i} \ln( v_{i}) d\hat{\nu} +  \int P_{v} \ln(\phi) d\hat{\nu}=$$
$$= - \int \sum_{i=0}^{d-1}  v_{i} \ln( v_{i}) d\hat{\nu} +  \int \sum_{i=0}^{d-1}
 v_{i} \ln( \phi(\tau_{i})) d\hat{\nu}=$$
$$= - \int \sum_{i=0}^{d-1}  v_{i} \ln( v_{i}) d\hat{\nu} +  \int \sum_{i=0}^{d-1}
 v_{i} \ln( \frac{\phi(\tau_{i}) h(\tau_{i})}{\rho h} \frac{\rho h}{h(\tau_{i})}) d\hat{\nu}=$$
$$= - \int \sum_{i=0}^{d-1}  v_{i} \ln( v_{i}) d\hat{\nu} +  \int \sum_{i=0}^{d-1}
 v_{i} \ln( u_{i}  \frac{\rho h}{h(\tau_{i})}) d\hat{\nu}=$$
$$=  \int -\sum_{i=0}^{d-1}  v_{i} \ln( v_{i}) + \sum_{i=0}^{d-1}  v_{i} \ln( u_{i} ) d\hat{\nu}
+ \int \sum_{i=0}^{d-1}  v_{i} \ln( \frac{\rho h}{h(\tau_{i})})
d\hat{\nu} \leq $$
$$= \ln(\rho) + \int \sum_{i=0}^{d-1}  v_{i} \ln(h) d\hat{\nu} - \int \sum_{i=0}^{d-1}
 v_{i} \ln(h(\tau_{i})) d\hat{\nu} =  $$
$$= \ln(\rho) + \int \ln(h) d\hat{\nu} - \int P_{v}\ln(h) d\hat{\nu} =  \ln(\rho)$$

The equality follows from the Lemma~\ref{PositOfEntropyAlt}
\end{proof}

From Theorem~\ref{PressureIsLogRoAlt} and
Lemma~\ref{PositOfEntropyAlt}, it follows that there exists
equilibrium states, more precisely, given a $\rho$-weighted
system, all holonomic liftings of the normalized probability, are
equilibrium states.

The Variational principle in the  formulation of the present section is stronger than
the formulated in the first part. The change in the definition
of entropy allow us to get a characterization of the equilibrium states
as holonomic liftings of the $P_{u}$-invariant probabilities of the
normalized transference operator. This point will become clear in the proof (of the "if, and only if," part) of the next theorem.

\begin{theorem}(Alternative Variational principle)\label{VarPrincipleAlt}
  Let us consider $\phi \in \mathbb{B}^{+}$ such that $([0,1], \tau_{i}, \phi)$
  is a $\rho$-weighted system, for $\rho = e^{p(\phi)} \geq 0$. Then, a holonomic
  probability $\hat{\nu}_{0}$ is a equilibrium state, if and only if,
  the projection  $\hat{\nu}_{0}$ by disintegration, is invariant  for
  the normalized dual operator of $P_{\phi}$( that is, $P_{u}^{*} (\pi_{*}\hat{\nu}_{0})=
  \pi_{*}\hat{\nu}_{0}$).
\end{theorem}
\begin{proof}
By Lemma~\ref{PositOfEntropyAlt}, it follows that: if
$\hat{\nu}_{0}$ is the  holonomic lifting of the normalized
probability $\pi_{*}\hat{\nu}_{0}$, then  $\hat{\nu}_{0}$ is an
equilibrium state.

The converse is also true. In fact, suppose that $\hat{\nu}_{0}$
is a equilibrium state, that is,
$$h(\hat{\nu}_{0}) + \int \ln(\phi) d\hat{\nu_{0}} = \ln(\rho)$$
Using the normalization we get,
$$- \int \sum_{i=0}^{d-1}  v_{i} \ln( v_{i}) d\hat{\nu_{0}} +
\int \ln(\phi) d\hat{\nu_{0}} = \ln(\rho)$$ where  $([0,1],
\tau_{i}, v_{i})$ is  the 1-weighted system associated to
holonomic disintegration of $\hat{\nu}_{0}$. From the invariance
of $P_{v}$ we have
$$- \int \sum_{i=0}^{d-1}  v_{i} \ln( v_{i}) d\hat{\nu_{0}} +
\int \sum_{i=0}^{d-1}  v_{i} \ln( \phi(\tau_{i})) d\hat{\nu_{0}} =
\ln(\rho)$$

Finally,  from the relations of normalization $P_{u}$ of
$P_{\phi}$
$$u_{i}(x)= \frac{\phi(\tau_{i}(x)) h(\tau_{i}(x))}{\rho h(x)}
\Leftrightarrow  \phi(\tau_{i}(x))= \frac{u_{i}(x) \rho
h(x)}{h(\tau_{i}(x))},$$ we get
$$- \int \sum_{i=0}^{d-1}  v_{i} \ln( v_{i}) d\hat{\nu_{0}} +
\int \sum_{i=0}^{d-1}  v_{i} \ln( \frac{u_{i}(x) \rho
h(x)}{h(\tau_{i}(x))}) d\hat{\nu_{0}} = \ln(\rho).$$

This is equivalent to
$$\int \sum_{i=0}^{d-1}  v_{i} \ln( \frac{u_{i}}{v_{i}}) d\hat{\nu_{0}} = 0.$$

From, $$\sum_{i=0}^{d-1}  v_{i} \ln( \frac{u_{i}}{v_{i}}) \leq
\ln(\sum_{i=0}^{d-1}  v_{i} \frac{u_{i}}{v_{i}}) =
\ln(\sum_{i=0}^{d-1}  u_{i})=\ln(1)=0 ,$$ it follows that $ u_{i}
= v_{i}, \quad \pi_{*}\hat{\nu}_{0}-a.e$.

As, $\pi_{*}\hat{\nu}_{0}$ is $P_{v}$-invariant, we get $P_{u}^{*}
(\pi_{*}\hat{\nu}_{0})= \pi_{*}\hat{\nu}_{0}$

\end{proof}

Finally, we point out that for $\phi$ fixed  one can consider a real parameter $\beta$
and the problem
$$p( \phi^{\beta})= \sup_{\hat{\nu} \in \mathcal{H}} \{ h(\hat{\nu}) +
    \,\beta\,  \int \ln(\phi) d\hat{\nu} \}.$$

For each value $\beta$,     denote by $\hat{\nu}_\beta$ a solution (therefore, normalized) of the above variational problem. Any subsequence  (weak limit) $\hat{\nu}_{\beta_n}\to \nu$ will determine a maximizing holonomic probability $\nu$ (in the sense, of maximizing $\sup_{\hat{\nu} \in \mathcal{H}} \{\,\int \ln(\phi) d\hat{\nu} \,\}$)  because the  entropy of any holonomic probability is bounded by $\ln d$. We refer the reader to \cite{GL} for properties on maximizing holonomic probabilities and we point out that  these results apply also for  the iterated setting as described above in the first two sections.

\vspace{0.5cm}

Instituto de Matem\'atica - UFRGS

Av. Bento Gon\c calves, 9500

Porto Alegre - 91500 - RS

Brasil
\vspace{0.2cm}

artur.lopes@ufrgs.br
\vspace{0.2cm}

oliveira.elismar@gmail.com

\vspace{0.5cm}

\end{document}